\documentclass[12pt]{article}
\usepackage{amssymb}
\usepackage[cp1251]{inputenc}
\usepackage[english]{babel}
\usepackage{amsmath}

\def\mapr#1{\smash{\mathop{\buildrel{#1}\over\longrightarrow}}}

\newlength{\myunit}
\newcounter{mycount}
\def\R{{\bf R}}
\def\F{{\bf F}}

\def\E{{\bf E}}
\def\Pr{{\bf Pr}}

\def\rank{\mbox{rank}}

\def\const{\mbox{\bf const }}

\def\inxx#1{}

\newtheorem{thm}{Theorem}

\newtheorem{df}{Definition}

\newtheorem{cor}{Corollary}
\newtheorem{rem}{Remark}

\newcommand\proof{{\bf Proof. }\nobreak\noindent}

\def\mathrm#1{\hbox{\rm #1 }}

\author{Anwar A. Irmatov}
\title{A Lower Bound of the Number of Threshold Functions in Terms of Combinatorial Flags on the Boolean Cube} 

\date{}
\begin{document}

\renewcommand\refname{\centering  \sc References }

\maketitle
\begin{abstract}
Let $ E=\{ (1, b_1, \ldots , b_n)\in {\R}^{n+1} \mid \; b_i= \pm 1 ,\; i=1, \ldots,  n \}$,  
$$E^{\times n}_{\ne 0} := \{ W=(w_{i_1}, \ldots , w_{i_n}) \mid w_{i_k}\in E, \, k=1, \ldots, n, \, dim \, span(w_{i_1}, \ldots , w_{i_n}) = n \},$$ 
and
$q^W_l := |span(w_{i_{n-l+1}}, \ldots , w_{i_n}) \cap E|.$
Then for any weights \\
$p=(p_1, \ldots, p_{2^n})$, $p_i\in \R$, $\sum_{i=1}^{2^n}{p_i} =1$ we have for the number of threshold functions $P(2,n)$ the following lower bound

$$P(2, n) \geq 2\sum_{W\in E^{\times n}_{\ne 0}}{\frac{1- p_{i_1} -p_{i_2} - \cdots - p_{i_{q_n^W}}}{q_n^W\cdot q_{n-1}^W\cdots q_1^W}},$$
\noindent and the right side of the inequality doesn't depend on the choice of $p$.
Here the indices used in the numerator correspond to vectors from $span(w_{i_1}, \ldots , w_{i_n})\cap E = \left\{w_{i_1}, \ldots, w_{i_n}, \ldots w_{i_{q_n^W}}\right\}$.

\bigskip

{\bf Keywords}\footnote{ This research has been partially supported by RFBR grant 18-01-00398 A}: threshold function, combinatorial flag, M\"obius function.
\end{abstract}


\section{Introduction.}

Let $ E_2 =\{\pm 1\}$ and $E^n_2:= \underbrace{E_2\times \cdots \times E_2}_{n}$.

\begin{df}\label{df1}
A function $f: E^n_2 \to E_2 $ is called a threshold function if there exist real numbers $\alpha_0, \alpha_1, \ldots , \alpha_n$
such that
$$ f(x_1, \ldots , x_n) =1 \mbox{\quad iff \quad} \alpha_1x_1+ \cdots + \alpha_n x_n + \alpha_0 \geq 0.$$
\end{df}

\noindent Denote by $P(2,n)$ the number of threshold functions.

Let us note that
$$f(x_1,\ldots , x_n) = \mbox{sign}  \langle \bar \alpha , (1, \bar x) \rangle,$$
where
\\
$(1, \bar x) := (1, x_1, \ldots, x_n) \in \R^{n+1} \mbox{\quad and \quad} \bar \alpha = (\alpha_0, \ldots , \alpha_n) \in \R^{n+1}.$
\\
\noindent This observation allows us to identify a threshold function with its $(n+1)$- weight vector $\bar \alpha$ which is a point in the dual space $(\R^{n+1})^{\ast} = \R^{n+1}$. It is shown in the paper {\rm \cite{Wi2}} that $P(2,n)$ can be expressed by the number $C(E)$ of disjoint chambers obtained as a compliment in $\R^{n+1}$ to the arrangement of $2^n$ hyperplanes all passing through the origin with the normal vectors from the set

\begin{equation}\label{eq1}
E=\{ (1, b_1, \ldots , b_n) \mid \; b_i \in E_2 ,\; i=1, \ldots,  n \}.
\end{equation}

The upper bound of the number $C(H)$ for any central arrangement of hyperplanes with a set $H$ of normal vectors was establisched by L.~Schl\"afli in {\rm \cite{Sch}}. For the case $H=E$ we have the following upper bound:

\begin{equation}\label{eq2}
P(2, n) = C(E) \leq 2 \sum_{i=0}^n {2^n-1 \choose i}.
\end{equation}

It should be noted that in the early 60s of the $20^{th}$ century the upper bound (\ref{eq2}) was obtained by several authors {\rm \cite{Cam}, \cite{Jos}, \cite{WhW}, \cite{Wi1}}. The detailed information of contribution of above mentioned authors can also be found in {\rm \cite{Cov}}.

One of the first lower bound of $P(2, n)$ was established by S.~Muroga in {\rm \cite{Mur}}:

\begin{equation}\label{eq3}
P(2, n) \geq 2^{0.33048 n^2}.
\end{equation}

S. Yajima and T. Ibaraki in {\rm \cite{YaI}} improved the order of the logarithm of the lower bound (\ref{eq3}) upto $n^2/2 $ :

\begin{equation}\label{eq4}
P(2, n) \geq 2^{n(n-1)/2 +8} \mbox{\quad for \quad } n \geq 6.
\end{equation}

Further significant improvements of the bound (\ref{eq4}) were obtained basing on the  paper {\rm \cite{Odl}} of A.~M.~Odlyzko. In the paper {\rm \cite{Zue}} it was noted that from the papers {\rm \cite{Odl}, \cite{Zas}} follows:

\begin{equation}\label{eq5}
C(E) =P(2, n) \geq 2^{n^2 - 10n^2/\ln n +O(n \ln n)} .
\end{equation}

Taking into account the upper bound (\ref{eq2}) and the inequality (\ref{eq5}) it is easy to see that

\begin{equation}\label{eq6}
\lim_{n \to \infty}{\frac{\log_2 P(2, n)}{n^2}} = 1.
\end{equation}

 In the paper {\rm \cite{Ir1}} a combination of an original geometric construction with the result from the paper {\rm \cite{Odl}} made it possible to improve the inequality (\ref{eq5}) upto:

\begin{equation}\label{eq7}
P(2, n) \geq 2^{n^2\left(1-\frac{7}{\ln n}\right)}\cdot  P\left(2, \left[\frac{7(n-1)}{\log_2 (n-1)}\right]\right) .
\end{equation}

The generalization of the inequality (\ref{eq7}) for the number of threshold $k$-logic functions was obtained in {\rm \cite{IrK}}. Asymptotics of logarithm of the number of polinomial threshold functions has been recently obtained in {\rm \cite{BV}}.

Let $H = \{ w_1, \ldots , w_T \}$ be a set of vectors in $\R^{n+1}$ . 
Let us denote by $H^{\times s}$ the set of ordered collections of different $s$  vectors, $s=1,\ldots, T$, from $H$ and let $[H]^{\times s}_{\ne 0}\subset H^{\times s}$ be a subset of linearly independent vectors
$$H^{\times s}_{\ne 0} := \{ (w_{i_1}, \ldots , w_{i_s}) \mid dim \, span(w_{i_1}, \ldots , w_{i_s}) = s \}.$$

For any $W= (w_{i_1}, \ldots , w_{i_n}) \in H^{\times n}_{\ne 0}$ let 
\begin{equation}\label{eq27}
q^W_l :=|L_l(W) \cap H| := |span(w_{i_{n-l+1}}, \ldots , w_{i_n}) \cap H|.
\end{equation}

\begin{df}\label{df2}
For any $W \in H^{\times n}$ the ordered set of numbers
$$W(H):=(q^W_n, q^W_{n-1}, \ldots, q^W_1)$$
is called a combinatorial flag of the ordered set $W$.

If $W \in [H]^{\times n}_{\ne 0}$, then $W(H)$ is called a full combinatorial flag of $W$.
\end{df}
For the sake of simplicity we will use the following notation:
\begin{equation}\label{A2}
W[H]:=q^W_n\cdot q^W_{n-1} \cdots q^W_1.
\end{equation}
The aim of this paper is to prove the following lower bound of the number of threshold functions. For any weights $p=(p_1, \ldots, p_{2^n})$, $p_i\in \R$, $\sum_{i=1}^{2^n}{p_i} =1$ on the set $E$ (see (\ref{eq1})) we have

\begin{equation}\label{eqA1}
P(2, n) \geq 2\sum_{W\in E^{\times n}_{\ne 0}}{\frac{1- p_{i_1} -p_{i_2} - \cdots - p_{i_{q_n^W}}}{W[E]}},
\end{equation}
and the right side of the inequality doesn't depend on the choice of $p$.
Here the indices used in the numerator correspond to vectors from $L_n(W)\cap E = \left\{w_{i_1}, \ldots, w_{i_n}, \ldots w_{i_{q_n^W}}\right\}$.

\section{Preliminaries.}

Let $H^{\perp}$ be a finite arrangement of hyperplanes all passing through the zero in $\R^{n+1}$ (central arrangement) and denote by $H = \{ w_1, \ldots , w_T \}$ the set of their normal vectors. We define a partially ordered set (poset) $L^H$ in the following way. By definition any subspace of $\R^{n+1}$ generated by some (possibly empty) subset of $H$ is an element of the poset $L^H$. An element $s \in L^H$ is less than an element $t \in L^H$ iff the subspace $t$ contains the subspace $s$. For any poset $P$ we can define a {\it simplicial complex} $\Delta_P$ in the following way. The set of vertices of $\Delta_P$ coincides with the set of elements $P$ and a set of vertices of $P$ defines a simplex of $\Delta_P$ iff this set forms a chain in $P$. Let us denote by $\Delta_{L^H}$ the simplicial complex of the poset
$$\left( 0_{L^H} , 1_{L^H}\right) := \{ z \in L^H \mid 0_{L^H} < z < 1_{L^H} \} ,$$
\noindent where $0_{L^H}$ and $1^{L^H}$ are the elements of the poset $L^H$ corresponding to the zero subspace of $\R^{n+1}$ and the whole space $\R^{n+1}$, respectively. It has been shown in {\rm \cite{Zas}} that the number $C(H)$ of $(n+1)$-dimensional regions into which $\R^{n+1}$ is divided by hyperplanes from the set $H^{\perp}$ can be found by the formula:
\begin{equation}\label{eq17}
C(H) = \sum_{t \in L^H}{\left|  \mu(0_{L^H} , t)\right|},
\end{equation}
\noindent where $\mu (s , t)$ is M\"obius function of the poset $L^H$.

M\"obius functions of partially ordered sets in the Zaslavsky's formula (\ref{eq17}) for calculation of the number of chambers $C(H)$ can be interpreted by tools of algebraic topology in the following way.
First, we introduce a simplicial compex $K^H$. The set of vertices of $K^H$ coincides with the set $H$. A subset $\{w_{i_1}, \ldots ,w_{i_s} \}$ of $H$
forms a simplex of $K^H$ iff
$$span(w_{i_1}, \ldots ,w_{i_s}) \ne \R^{n+1}.$$

Taking into account the results of the papers {\rm \cite{Fol}, \cite{Hal}} it is possible to show (see {\rm \cite{Ir3}}) that the absolute value of the M\"obius function $|\mu (0_{L^H}, u)|$  is equal to the dimension of the reduced homology group of the complex $K^{H\cap u}$ with coefficients in an arbitrary field $\F$:
\begin{equation}\label{eq18}
|\mu (0_{L^H}, u)| = \rank \, \tilde H_{dim\, u -2}\left( K^{H\cap u} ; \F\right).
\end{equation}
\noindent Here the set $H\cap u$ consists of all vectors $H$ belonging to the subspace $u \subset \R^{n+1}$ and is considered as a subset of $\R^{dim\, u}:= u.$

It follows from the definition of M\"obius function that
\begin{equation}\label{eq19}
\sum_{0_{L^H} \leq u < 1_{L^H}}{|\mu (0_{L^H} , u)|} \geq \left| -\sum_{0_{L^H} \leq u < 1_{L^H}}{\mu (0_{L^H} , u)}\right| = |\mu (0_{L^H} , 1_{L^H})|.
\end{equation}
\noindent Hence,

\begin{equation}\label{eq20}
C(H) = |\mu (0_{L^H} , 1_{L^H})| + \sum_{0_{L^H} \leq u < 1_{L^H}}{|\mu (0_{L^H} , u)|} \geq 2|\mu (0_{L^H}, 1_{L^H})|.
\end{equation}

\noindent From (\ref{eq18}) and (\ref{eq20}) we have:
\begin{equation}\label{eq21}
C(H) \geq 2\,\rank \, H_{n-1}\left( K^H ; \F\right).
\end{equation}

\noindent As a consequence of (\ref{eq21}) for the case $H=E$, we have:
\begin{equation}\label{eq22}
P(2 , n) = C(E) \geq 2\, \rank\, H_{n-1}\left( K^{E} ; \F\right).
\end{equation}

Let us fix an arbitrary order $\pi : [T] \to H$ on the set $H$.  We denote by $\Lambda^H$ the number of all collections of different vectors $(w_{i_1}, \ldots , w_{i_n})$, $2\leq i_1, \ldots , i_n \leq T$ such that for any $l, \; 1\leq l \leq n,$ the vector $w_{i_l}$ is the minimal vector in the order $\pi$ among all vectors from the set $H\cap span (w_{i_l}, \ldots , w_{i_n}).$ It follows from the definition of $\Lambda^H$ that 
$$dim\,span(w_{i_l}, \ldots , w_{i_n}) = n-l+1,$$
\noindent i.e. the vectors $w_{i_1}, \ldots , w_{i_n}$ are linearly independent.

The theorem~7 of {\rm \cite{Ir3}} is also true for any set $H$ generating $\R^{n+1}$ and here we rewright it as follows:
\begin{equation}\label{eq23}
\rank\, H_{n-1}\left( K^H ; \F\right) = \Lambda^H.
\end{equation}

The description (\ref{eq23}) of the number $\Lambda^H$ as rank of the homology group shows us that it doesn't depend on the order $\pi$ on $H$ and we can fix the order $\pi$  as $\pi (i):= w_i, \, i\in[T]$. In that case $w_i <_{\pi} w_j \iff i<j.$ We denote by $\Gamma$ 
the set of all orders on the set $H$. Then any order on $H$ can be defined as composition
$$[T] \mapr{\gamma} [T] \mapr {\pi} H$$
\noindent of a permutation $\gamma: [T] \to [T]$ with $\pi$, and 
$$w_i <_{\gamma} w_j \iff (\pi \gamma)^{-1}(w_i) < (\pi \gamma)^{-1}(w_j) \iff \gamma^{-1}(i) < \gamma^{-1}(j).$$

Thus, $\Gamma$ can be identified with the symmetric group $Sym([T])$ and any permutation $\sigma : [T] \to [T]$ defines a basis of the homology group $H_{n-1} \left(K^H ; \F\right)$ considering as a vector space over the field $\F$ as a subset of collections of $n$ vectors from $H$
$$B^{\sigma} \subset H^{\times s}_{\ne 0} $$
\noindent obeying the following conditions. Let $W^{\sigma} = (w_{\sigma(i_1)}, \ldots , w_{\sigma (i_n)})$ and $L(W^{\sigma})$ be a flag of subspaces in $\R^{n+1}$
$$L_n(W^{\sigma}) \supset L_{n-1}(W^{\sigma}) \supset \ldots \supset L_l(W^{\sigma}) \supset \ldots \supset L_1(W^{\sigma}),$$
\noindent where
\begin{equation}\label{eq24}
L_l(W^{\sigma}) := span( w_{\sigma(i_{n-l+1})}, \ldots , w_{\sigma (i_n)}) \quad l=1, \dots, n.
\end{equation}
Then $W^{\sigma} \in B^{\sigma}$ iff 

\begin{equation}\label{eqA25}
 i_1< i_2 < \ldots < i_n
\end{equation}

\begin{equation}\label{eq25}
w_{\sigma(1)} \notin L_n(W^{\sigma})
\end{equation}
\noindent and
\begin{equation}\label{eq26}
\min \{ i\in [T] \mid w_{\sigma(i)} \in H\cap L_l(W^{\sigma})\} = i_{n-l+1}, \quad l=1, \ldots, n.
\end{equation}

\bigskip

\section{A lower bound of threshold functions.}

\begin{thm}\label{th1}
For any probability distribution $p=(p_1, \ldots , p_T)$ on the set $H$ the following equality is true: 
$$\rank \, H_{n-1}\left( K^H ; \F \right) = \sum_{W\in H^{\times n}_{\ne 0}}{\frac{1- p_{i_1} -p_{i_2} - \cdots - p_{i_{q_n^W}}}{W[H]}}.$$
Here the indices used in the numerator correspond to vectors from $L_n(W)\cap H = \left\{w_{i_1}, \ldots, w_{i_n}, \ldots w_{i_{q_n^W}}\right\}.$
\end{thm}

\proof  We define a probability distribution $\tilde p$ on the set $\Gamma \cong Sym([T])$ by the formula:
$$\tilde p(\gamma) =p_{\gamma(1)}\frac{1}{(T-1)!} , \quad \gamma \in Sym([T]).$$

For any collection $W=(w_{i_1}, \ldots , w_{i_n}) \in H^{\times n}_{\ne 0}$ we define a random function $I_W : \Gamma \to \R$ by the formula:

$$I_W(\gamma) = \left\{\begin{array}{ll}
1, & \mbox{if\;} \gamma^{-1}(i_1) < \gamma^{-1}(i_2) < \cdots < \gamma^{-1}(i_n) \\
    & \mbox{and the conditions (\ref{eq25}), (\ref{eq26}) hold for the collection} \\
    &  W^{\gamma}:=(w_{\gamma (\gamma^{-1}(i_1))}, \ldots, w_{\gamma(\gamma^{-1}(i_n))}) = (w_{i_1}, \ldots , w_{i_n}); \\
    & \mbox{in the order generated by \;} \gamma \\
0, & \mbox{in all other cases.}
\end{array}\right. $$

Let
$$I := \sum_{W\in H^{\times n}_{\ne 0}}{I_W} : \Gamma \to \R.$$

Then for any $\gamma \in \Gamma$
$$I(\gamma) = \const = \rank \, H_{n-1}\left(K^H ; \F\right).$$
Hence the expectation of $I$ is equal to the rank of the homology group:
\begin{equation}\label{eq28}
\E [I] = \rank\, H_{n-1}\left(K^H ; \F\right).
\end{equation}
Additivity of expectation reduces the problem of calculation $\E[I]$ to counting the probability $\Pr (I_W =1)$:
\begin{equation}\label{eq29}
\E[I] = \sum_{W\in H^{\times n}_{\ne 0}}{\E[I_W]} = \sum_{W\in H^{\times n}_{\ne 0}}{\Pr (I_W =1)}.
\end{equation}

Further we calculate the number of permutations $\gamma$ such that $I_W(\gamma)=1$. Since $w_{\gamma(1)}\notin L_n(W)$ then $q^W_n$ vectors from 
$L_n(W)\cap H$ can be located in any places except the first one, i.e. $\gamma^{-1}(j) \ne 1$ for any $j\in [T]$ such that $w_j \in L_n(W)\cap H.$ The arrangement of the remaining vectors from $H\backslash \{w_{\gamma (1)}\cup \{L_n(W)\cap H\}$ does not affect the fulfillment of the conditions (\ref{eq25}) and (\ref{eq26}). By the condition (\ref{eq26}) the vector $w_{i_1}$ has to be located in the first place from the chosen $q^W_n$ positions for arrangement of the set $L_n(W)\cap H,$ while $q_{n-1}^W$ vectors from $L_{n-1}(W)\cap H$ can be located in any of the remaining $q_n^W-1$ places. The arrangement of the vectors from $\{L_n(W)\cap H\} \backslash \{w_{i_1}\cup \{L_{n-1}(W)\cap H\}\}$ in $q_n^W-q_{n-1}^W-1$ places left after choosing $q_{n-1}^W+1$ places for arrangement of the set $L_{n-1}(W)\cap H$ and $w_{i_1}$ doesn't affect the fulfillment of the conditions (\ref{eq25}) and (\ref{eq26}). Continuing the same way, we get that in the first place from $q^W_l$ positions arranged for the vectors from $L_l(W)\cap H$ has to be located the vector $w_{i_{n-l+1}}$, while $q^W_{l-1}$ vectors from $L_{l-1}(W)\cap H$ can be located in any of the remained $q^W_l-1$ places and the positions of the vectors $\{L_l(W)\cap H\} \backslash \{w_{i_{n-l+1}}\cup \{L_{l-1}(W)\cap H\}\}$ in $q_l^W-q_{l-1}^W-1$  places left after choosing $q_{l-1}^W+1$ places for arrangement of the set $L_{l-1}(W)\cap H$ and $w_{i_{n-l+1}}$ doesn't affect the fulfillment of the conditions (\ref{eq25}) and (\ref{eq26}).

Denote by $N(\gamma(1)=i)$ the number of permutations $\gamma$ with fixed value $\gamma(1)=i$ such that $w_i \notin L_n(W)$. Then
$$\begin{array}{l} 
N(\gamma(1)=i) = 
{T-1 \choose q_n^W} (T-1-q_n^W)!\cdot {q_n^W-1 \choose q_{n-1}^W}(q_n^W-q_{n-1}^W-1)!\cdots \\
\\
\cdot {q_l^W-1 \choose q_{l-1}^W}(q_l^W-q_{l-1}^W-1)!\cdots {q_2^W-1 \choose q_1^W}(q_2^W-q_1^W-1)!(q_1^W-1)!=  \\
\\
= \frac{(T-1)!}{q_n^W!}\cdot \frac{(q_n^W-1)!}{q_{n-1}^W!}\cdots \frac{(q_l^W-1)!}{q_{l-1}^W!}\cdots \frac{(q_2^W-1)!}{q_1^W!}(q_1^W-1)!  
= \frac{(T-1)!}{q_n^W q_{n-1}^W \cdots q_2^W\cdot q_1^W} = \\
\\
=\frac{(T-1)!}{W[H]}
\end{array}$$
\noindent Then, we have
$$\Pr (I_W=1) = \sum_{i\in [T]\, s.t. \, w_i\notin L_n(W)}{p_i\frac{1}{(T-1)!}\frac{(T-1)!}{W[H]}} = $$
$$= \frac{1-p_{i_1} -\cdots p_{i_{q_n^W}}}{W[H]},$$
\noindent where $ L_n(W)\cap H= \{w_{i_1},\ldots, w_{i_n}, \ldots, w_{i_{q_n^W}}\}.$

\begin{flushright} {\sc Q.E.D.} \end{flushright}

\begin{rem}
Since on the right side of equation of Theorem \ref{th1} is a polynomial of degree 1 then the Theorem \ref{th1} is true for any $p_i \in \R$, $i=1,\ldots, T$, such that $ \sum_{i=1}^T{p_i}=1$.
\end{rem}

\begin{cor}
For any weights $p=(p_1, \ldots, p_{2^n})$, $p_i\in \R$, $\sum_{i=1}^{2^n}{p_i} =1$ on the set $E$ we have

$$P(2, n) \geq 2\sum_{W\in E^{\times n}_{\ne 0}}{\frac{1- p_{i_1} -p_{i_2} - \cdots - p_{i_{q_n^W}}}{W[E]}},$$

\noindent and the right side of the inequality doesn't depend on the choice of $p$.
Here the indices used in the numerator correspond to vectors from $L_n(W)\cap E = \left\{w_{i_1}, \ldots, w_{i_n}, \ldots w_{i_{q_n^W}}\right\}$.
\end{cor}
\proof The enequality follows from (\ref{eq22}) and Theorem \ref{th1}.

\bigskip

\noindent Faculty of Mechanics and Mathematics, Lomonosov MSU, Moscow, Russian Federation

\bigskip

{\it E-mail address}: irmatov@intsys.msu.ru

\end{document}